\documentclass[draftclsnofoot,onecolumn]{IEEEtran}

\usepackage{verbatim}
\usepackage[inline]{enumitem1}
\usepackage{setspace,mathtools,graphicx,amsfonts}

\graphicspath{{images/}}

\newtheorem{mytheorem}{Theorem}
\newtheorem{myexample}{Example}
\newtheorem{mydef}{Definition}
\newtheorem{mylemma}{Lemma}
\newtheorem{mycorollary}{Corollary}

\begin{document}

\title{Fa\`a di Bruno's formula for G\^ateaux differentials and interacting stochastic population processes}

\author{Daniel~E.~Clark and~Jeremie~Houssineau}

\maketitle

\begin{abstract}
The problem of estimating interacting systems of multiple objects is important to a number of different fields of mathematics, physics, and engineering. 
Drawing from a range of disciplines, including statistical physics, variational calculus, point process theory, and statistical sensor fusion, we develop a unified probabilistic framework for modelling systems of this nature.
In order to do this, we derive a new result in variational calculus, Fa\`a di Bruno's formula for G\^ateaux differentials.
Using this result, we derive the Chapman-Kolmogorov equation and Bayes' rule for stochastic population processes with interactions and hierarchies. 
We illustrate the general approach through case studies in multi-target tracking, branching processes and renormalization.
\end{abstract}

\section*{Introduction}

In many science and engineering applications, researchers are interested in estimating the state of dynamical systems that have uncertainty in their dynamics and population.
Stochastic filtering methods are of particular importance, which are modelled for discrete-time systems using the Chapman-Kolmogorov equation \cite{chapman,kolmogorov} and Bayes' rule~\cite{Bayes}.
The extension of stochastic filtering to systems of multiple object systems is a recent development, developed by Mahler \cite{Mahlersteinwinter} as a means of identifying and tracking an unknown number of objects in aerospace applications. 
This was based on the Chapman-Kolmogorov equation for stochastic population processes proposed by Moyal \cite{Moyal62} and the derivation of Bayes' rule for point processes by Mahler \cite{Mahlermaths}.
Since then, this approach has been widely used to solve multi-object estimation problems \cite{VoMa1,VoVoCantoni,VoClarkRistic,RistikClarkVoVo}.
In this paper, we develop a class of  processes for modelling systems of multiple objects to account for dependencies between groupings of objects inspired from the idea proposed by Ursell \cite{ursell},
and used to determine:
\begin{enumerate}[label=(\roman{*})]
 \item A means of working with point processes and stochastic population processes when the associated probability measure is not assumed symmetric.
 \item A generalisation of  branching processes~\cite{harris,MoyalAP,MoyalRS} that models correlations.
 \item A definition of the problem known as renormalization in percolation theory~\cite{Stauffer} in terms of stochastic population processes.
 \item A Bayesian solution to the problem of multi-object estimation~\cite{Mahlermaths} with applications. 
\end{enumerate}

\subsection*{Stochastic population processes and multi-object estimation}

The mathematical foundation for stochastic population processes~\cite{Moyal62} is based on the theory of generating functionals~\cite{Volterra}.
The probability generating functional~\cite{Moyal62,sb}, p.g.fl., can be viewed as a generalisation of the more commonly known probability generating function~\cite{Galton}, p.g.f., for modelling population processes.
In a similar manner to the p.g.f., the probability density and factorial moment densities of multi-object processes can be found from the p.g.fl.\ by differentiation, with G\^{a}teaux differentials~\cite{Gateaux}.
The probability generating functional is well known within the point process literature~\cite{DaleyVere-Jones,RCoxIsham} yet the use of this concept is often restricted to the construction or specification of basic models, since it can lead to quite cumbersome formulae, and the  generality of the generating functional approach for stochastic modelling is still yet to be fully exploited.

Mahler \cite{Mahlermaths,MahlerCPHD} was the first to make extensive use of the p.g.fl.\ for stochastic modelling of multi-object  estimation problems, within the framework of Finite Set Statistics (FISST), which was developed as a unified approach to multi-sensor multi-target data fusion~\cite{GoodmanMahlerNguyen,MahlerBook}.
This approach is now regularly used to derive practical multi-object estimation and sensor fusion algorithms, which are usually based on propagating an approximation of the first-order factorial moment density~\cite{Mahlersteinwinter}.

To detemine the multi-object analogue of Bayes' rule, Mahler \cite{Mahlermaths} proposed the use of functional derivatives of the probability generating functional.
This approach often involves finding the parameterised form of the updated process, and proving its correctness by induction.
The process can be quite involved and needs to be applied for each model~\cite{Mahlermaths,MahlerCPHD,swainspie10,Mahlerfusion2009,MahlerSPIE_09_3}.
The construction of the models often involves composition of basic models, whose derivatives are easy to find, yet when composed, their form for higher derivatives becomes more unclear.
In this paper we circumvent this problem by introducing a new tool in variational calculus, Fa\`a di Bruno's formula for G\^ateaux differentials on topological vector spaces.

\subsection*{Fa\`a di Bruno's formula}

Mathematicians have investigated formulae for expressing 
higher-order
derivatives of composite functions
in terms of derivatives of their factor functions
for over 200 years.
				These formulae are often attributed to  Fa\`{a} di Bruno~\cite{faa1,faa2},
				though  Craik~\cite{Craik}
				recently highlighted a number of researchers preceding
				his works, the earliest of which is thought to be 
				by Arbogast~\cite{arbo}.
				Despite the fact that the idea of expressing these formulae
				in terms of derivatives of the factor
				functions is not new, a number of recent works have appeared 
				on this topic,
				including those by Hardy~\cite{hardy} and Ma~\cite{ma}
				on partial derivatives,
				and an alternative approach was presented by Huang {\it et al.}
				for Fr\'{e}chet derivatives~\cite{huang}.

This paper presents a simple proof of the
formula to express the higher-order G\^{a}teaux differentials
of composite functions 
in terms of differentials of their factor functions.
This generalises the formula commonly attributed to
Fa\`{a} di Bruno 
to functions in locally convex topological vector spaces.
The result presented has  a simpler
form than previous formulae
and is described in terms of a sum over partitions of the increments.

\subsection*{Summary of results}

The paper is multi-disciplinary in nature and presents results in different domains.
In the next section, the general form of Fa\`a di Bruno's formula is derived. 
In Section \ref{sec:stoPopProc} the probability generating functional (p.g.fl.) is modified to account for non-symmetric point processes.
Interacting population processes are defined in Section \ref{sec:hierInterPopProc} through the application of Fa\`a di Bruno's formula on composite functionals in order to model a broad class of problems, including hierarchical population processes and processes with interactions.
Section \ref{sec:ChapmanKolmogorov} applies the interacting population processes as a means of unifying the problems of branching processes and renormalization as instances of the Chapman-Kolmogorov equation.
Section \ref{sec:BayesRule} presents the Bayesian filter and smoother in a general framework with a functional form that can account for interactions and hierarchies.
This is illustrated through the application to general Poisson and independent processes.
The paper concludes in Section \ref{sec:conclusion}.

\section{G\^ateaux differentials }

Since we require the use of differentials of functionals to determine the main results in the paper, we describe G\^ateaux differentials and their higher-order variations.
Moreover, we derive the most general form of Fa\`a di Bruno's formula.

The following two definitions describe the G\^ateaux differential and the $n^{th}$-order differential.
\begin{mydef}[G\^ateaux differential]
Let $X$ and $Y$ be locally convex topological vector spaces, let $\Omega$ be an open subset of $X$ and let $f:\Omega\rightarrow Y$. The G\^ateaux differential at $a\in \Omega$ in the direction $x\in X$ is
\begin{equation}
\delta f(a;x) = \lim_{\epsilon\to 0} \frac{1}{\epsilon}\left( f(a+\epsilon x) - f(a) \right)
\end{equation}
when the limit exists. If $\delta f(a;x)$ exists for all $x\in X$ then $f$ is Gateaux differentiable at $a$. The G\^ateaux differential is homogeneous of degree one in $x$, so that for all real numbers $\theta$, $\delta f(a;\theta x)=\theta\delta f(a;x)$.

\end{mydef}

\begin{mydef}[$n^{th}$-order G\^ateaux differential]
The {\it $n^{th}$-order variation} $\delta^nf(x;\eta_1,\ldots,\eta_n)$ of $f(x)$ in directions $\eta_1,\ldots,\eta_n \in X$ is defined recursively with
\begin{equation*}
\delta^n f\left(x;\eta_1,\ldots,\eta_n\right) =
\delta\left(\delta^{n-1} f\left(x;\eta_1,\ldots,\eta_{n-1}\right); \eta_n\right).
\end{equation*}
\end{mydef}

\subsection{Chain differentials}

Due to the lack of continuity properties of the G\^ateaux differential, we require further constraints in order to derive a chain rule.
Bernhard \cite{bernhard} proposed a new form of G\^ateaux differential defined with sequences, which he called the chain differential, that is not as restrictive as the Fr\'echet derivative though it is still possible to find a chain rule.

\begin{mydef}[Chain differential]
The function $f:X\rightarrow Y$, where $X$ and $Y$ are topological vector spaces, has a {\it chain differential} $\delta f(x;\eta)$ at $x$ in the direction $\eta$ if, for any sequence $\eta_n\rightarrow\eta\in X$, and any sequence of real numbers $\theta_n\rightarrow 0$, it holds that
\begin{equation*}
\delta f(x;\eta) = \lim_{n\rightarrow \infty} \dfrac{1}{\theta_n} \left( f(x+\theta_n\eta_n)-f(x) \right).
\end{equation*}
\end{mydef}

\begin{mylemma}[Chain rule, from \cite{bernhard}, Theorem 1]
\label{lemma:chainRule}
Let $X$, $Y$ and $Z$ be topological vector spaces, $f : Y \rightarrow Z$ , $g : X\rightarrow Y$ and $g$ and $f$ have chain differentials at $x$ in the direction $\eta$ and at $g(x)$ in the direction $\delta g(x;\eta)$ respectively. Let $h = f \circ g$, then $h$ has a chain differential at $x$ in the direction $\eta$, given by the chain rule
\begin{equation*}
\delta h(x;\eta) = \delta f(g(x); \delta g(x;\eta)).
\end{equation*}
\end{mylemma}

\begin{mydef}[Partial chain differential]
Let $\{X_i\}_{i=1:n}$ and $Y$ be topological vector spaces. The function $f:X_1\times\ldots\times X_n \rightarrow Y$ has a {\it partial chain differential} w.r.t. the $i^{th}$ variable $\delta_i f(x_1,\ldots,x_n;\eta)$ at $(x_1,\ldots,x_n)$ in the direction $\eta$ if, for any sequence $\eta_m\rightarrow\eta\in X$, and any sequence of real numbers $\theta_m\rightarrow 0$, it holds that
\begin{equation*}
\delta_i f(x_1,\ldots,x_n;\eta) = \lim_{m\rightarrow \infty} \dfrac{1}{\theta_m} \left( f(x_1,\ldots,x_i+\theta_m\eta_m,\ldots,x_n)-f(x_1,\ldots,x_n) \right).
\end{equation*}
\end{mydef}

{\begin{mytheorem}[Total chain differential]
\label{thm:totChainDiff}
Let $\{X_i\}_{i=1:n}$ and $Y$ be topological vector spaces. The function $f:X_1\times\ldots\times X_n \rightarrow Y$ has a {\it total chain differential} $\delta f$ at $(x_1,\ldots,x_n)$ if all the partial chain differentials exist in a neighbourhood $\Omega \subseteq X_1\times\ldots\times X_n$ of $(x_1,\ldots,x_n)$ and in any direction and if $\delta_i f$ is continuous over $\Omega\times X_i$. Then for $\underline{\eta} \in X_1\times\ldots\times X_n$ such that $\underline{\eta} = (\eta_1,\ldots,\eta_n)$,
\begin{equation*}
\delta f(x_1,\ldots,x_n;\underline{\eta}) = \sum_{i=1}^n \delta_i f(x_1,\ldots,x_n;\eta_i).
\end{equation*}

If $X_i = X$, $1\leq i \leq n$, and $\underline{\eta} = (\eta,\ldots,\eta)$, we write $\delta f(x_{1:n};\underline{\eta}) = \delta f(x_{1:n};\eta)$.
\end{mytheorem}}

\begin{IEEEproof}
The result is proved in the case $n=2$ from which the general case can be straightforwardly deduced:
\begin{align}
\label{eq:defTotalProof1}
\delta f\left(x,y; (\eta,\xi)\right) & = \lim_{r\rightarrow \infty} \theta^{-1}_r\left[ f(x+\theta_r \eta_r,y + \theta_r \xi_r) - f(x,y) \right]\\
\notag
& = \lim_{r\rightarrow \infty} (\theta^{-1}_r\left[f(x+\theta_r \eta_r,y + \theta_r \xi_r) - f(x+\theta_r \eta_r,y)\right] + \theta^{-1}_r\left[f(x+\theta_r \eta_r,y) - f(x,y)\right] ).
\end{align}

Given $\theta_r \neq 0$, define $h:\mathbb{R} \rightarrow \mathbb{R}$ as $h(t) = f(x+\theta_r \eta_r,y + t\xi_r)$. From the mean value theorem, there exists $c_y \in [0,\theta_r]$ such that $h(\theta_r) - h(0) = \delta h(c_y;\theta_r)$.
Using Lemma \ref{lemma:chainRule}, this last equation is equivalent to
\begin{equation*}
\theta^{-1}_r\left[f\left(x+\theta_r \eta_r,y + \theta_r\xi_r\right) - f\left(x+\theta_r \eta_r,y\right)\right] = \delta_2 f(x+\theta_r \eta_r,y + c_y \xi_r;\xi_r),
\end{equation*}
similarly for the partial chain differential w.r.t. the first variable, there exists $c_x \in [0,\theta_r]$ such that
\begin{equation*}
\quad \theta^{-1}_r\left[f(x+\theta_r \eta_r,y) - f(x,y)\right] = \delta_1 f(x+c_x \eta_r,y;\eta_r).
\end{equation*}
Because of the continuity of $\delta_1 f$ and $\delta_2 f$, (\ref{eq:defTotalProof1}) becomes
\begin{equation*}
\delta f\left(x,y; (\eta,\xi)\right) = \delta_1 f\left(x,y; \eta\right) + \delta_2 f\left(x,y; \xi\right)
\end{equation*}
which is equivalent to the Proposition 3 in \cite{bernhard}.
\end{IEEEproof}

\subsection{Fa\`a di Bruno's formula}

The next result generalises the chain rule for higher-order variations.
This is analogous to the formula of Fa\`a Di Bruno \cite{faa1} for G\^ateaux differentiation.
The important distinction between the chain rule for G\^ateaux differentials and other forms of differentiation is that the resulting formula does not factorise as a product.

\begin{mytheorem}[General chain rule]
\label{thm:chainRule}
Let $X$, $Y$ and $Z$ be topological vector spaces.
Assuming that and $g:X\rightarrow Y$ has higher order chain differentials in any number of directions in the set $\{\eta_1, \ldots, \eta_n\}\in X^n$ and that $f:Y\rightarrow Z$ has higher order chain differentials in any number of directions in the set $\{\delta^m g(x;S_m)\}_{m=1:n}$, $S_m \subseteq \{ \eta_1, \ldots, \eta_n \}$.
Assuming additionally that for all $1\leq m \leq n$, $\delta^m f(y;\xi_1,\ldots,\xi_m)$ is continuous on an open set $\Omega \subseteq Y^{m+1}$ and linear w.r.t. the directions $\xi_1,\ldots,\xi_m$, the $n^{th}$-order variation of composition $f\circ g$ in directions $\eta_1, \ldots, \eta_n$ at point $x\in X$ is given by
\begin{equation*}
\delta^n (f\circ g)(x; \eta_1,\ldots,\eta_n ) = \sum_{\pi\in \Pi(\eta_{1:n})} \delta^{\dot{\pi}}f \left( g(x); \xi_{\pi_1}(x),\ldots,\xi_{\pi_{\dot{\pi}}}(x) \right),
\end{equation*}
where $\xi_{\omega}(x) = \delta^{\dot{\omega}} g \left( x; \omega_1,\ldots,\omega_{\dot{\omega}} \right)$ is the $\dot{\omega}^{th}$-order chain differential of $g$ in directions $\{\omega_1,\ldots,\omega_{\dot{\omega}}\} \subseteq \{ \eta_1, \ldots, \eta_n \} $. $\Pi(\eta_{1:n})$ represents the set of partitions of the set $\{\eta_1,\ldots,\eta_n\}$ and $\dot{\pi}$ denotes the cardinality of the set $\pi$.
\end{mytheorem}

\begin{IEEEproof}
Lemma \ref{lemma:chainRule} gives the base case $n=1$.
For the induction step, we apply the differential operator to the case $n$ to give the case $n+1$ and show that it involves a summation over partitions of elements
$\eta_1,\ldots,\eta_{n+1}$
\begin{equation*}
\delta^{n+1} (f\circ g)(x; \eta_1,\ldots,\eta_{n+1}) = \sum_{\pi\in \Pi(\eta_{1:n})} \delta \left( \delta^{\dot{\pi}}f \left( g(x); \xi_{\pi_1}(x),\ldots,\xi_{\pi_{\dot{\pi}}}(x) \right); {\eta_{n+1}} \right).
\end{equation*}

Let $\phi^g_m$ and $\phi^{\xi_i}_m$, $1\leq i \leq n$ be defined as $\phi^h_m = \theta^{-1}_m \left( h(x + \theta_m \eta_{n+1}) - h(x) \right) $ so that $\phi^h_m\rightarrow \delta h\left(x; \eta_{n+1}\right)$ and $h(x + \theta_m \eta_{n+1}) = h(x) + \theta_m \phi^h_m (x)$. Let the multivariate function $F:Y^k \rightarrow Z$ be $F(h_0,\ldots,h_k) = \delta^k f \left( h_0; h_1,\ldots,h_k \right), \quad k \geq 0$.
Using Theorem \ref{thm:totChainDiff} and using the linearity of $\delta^k f$ w.r.t. to the directions $h_1,\ldots,h_k$, we write
\begin{align*}
\delta^{n+1} (f\circ g)(x;\eta_1,\ldots,\eta_{n+1}) & = \lim_{m\rightarrow\infty} \sum_{\pi\in \Pi(\eta_{1:n})}\delta F\left( g,\xi_{\pi_1},\ldots,\xi_{\pi_{\dot{\pi}}}; \phi^g_m,\phi^{\xi_{\pi_1}}_m,\ldots,\phi^{\xi_{\pi_{\dot{\pi}}}}_m\right)\\
& = \left.\sum_{\pi\in \Pi(\eta_{1:n})} \right[\delta^{\dot{\pi}+1}f \left( g; \xi_{\pi_1},\ldots,\xi_{\pi_{\dot{\pi}}},\delta g \right) + \left.\sum_{\omega\in\pi} \delta^{\dot{\pi}}f \left( g; \xi_{\pi_1},\ldots,\delta{\xi}_{\omega},\dots,\xi_{\pi_{\dot{\pi}}} \right)\right].
\end{align*}
where the argument of $g$ and $\xi_i$ has been omitted when there is no ambiguity.
The result above can be viewed as a means of generating all partitions of $n+1$ elements from all partitions of $n$ elements:
The first term takes the variation with increment $\delta g\left(x;\eta_{n+1}\right)$, and each term in the second summation takes the variation of each increment $\xi(x)$ in the direction $\eta_{n+1}$, i.e. $\delta{\xi}(x;\eta_{n+1})$.
\end{IEEEproof}

It is worth highlighting the structure of the result above.
In other forms of chain rule, Fa\`a di Bruno's formula is a sum over partitions of products.
However, in the general form above for G\^{a}teaux differentials, the outer functional has variations in directions that themselves are differentials of the inner functional.

\subsection{Useful instances of chain differentials}
The following two lemmas will be used in order to determine the probability and factorial moment measures of stochastic population processes.

\begin{mylemma}[Chain differential of a composed linear function]
\label{lemma:expectDiffInvers}
Let $F:Y\rightarrow Z$ be a linear function and $h:X\rightarrow Y$ be any function. The chain differential of an $F \circ h$ at point $x\in X$ in directions $\eta_i \in X$, $1\leq i \leq k$, satisfies the following equality
\begin{equation}
\delta^k \left(\left(F \circ h\right)(x);\eta_1,\ldots,\eta_k \right) = F\left( \delta^k h\left(x;\eta_1,\ldots,\eta_k \right) \right), \quad \forall k\in\mathbb{N}.
\end{equation}
\end{mylemma}

\begin{IEEEproof}
The case $k=0$ is obvious. Assuming the result is valid at order $k-1$, the $k^{th}$-order chain differential is $\delta (F( \delta^{k-1} h(x;\eta_1,\ldots,\eta_{k-1}));\eta_k)$.

We set $\bar{h}(x) = \delta^{k-1} h(x;\eta_1,\ldots,\eta_{k-1})$. From the chain rule (Lemma \ref{lemma:chainRule}), we find $\delta(F(\bar{h}(x));\eta_k) = \delta F(\bar{h}(x);\delta_1 \bar{h}( x;\eta_k))$ where
\begin{align*}
\delta F \left(\bar{h}\left(x\right);\delta_1 \bar{h}\left( x;\eta_k\right)\right) & = \lim_{m\rightarrow \infty} \theta^{-1}_m \left( F\left( \bar{h}(x) + \theta_m\delta \bar{h}\left(x;\eta_{k,m} \right) \right) - F\left( \bar{h}(x)\right) \right) \\
& = \lim_{m\rightarrow \infty} F\left(\delta \bar{h}\left(x;\eta_{k,m} \right) \right)\\
& = \lim_{m\rightarrow \infty} F\left(\delta \left( \delta^{k-1} h\left(x;\eta_1,\ldots,\eta_{k-1}\right) ;\eta_{k,m} \right) \right),
\end{align*}
which gives the desired result at order $k$ since $\displaystyle\lim_{m\rightarrow \infty}\eta_{k,m} = \eta_k$.
\end{IEEEproof}

\begin{mylemma}[Chain differentials of a multi-linear function]
\label{lemma:multiDiff}
Let $F:X_1\times\ldots\times X_n\rightarrow Y$ be a multi-linear function. The $k^{th}$-order partial chain differentials of $F$ at point $x=(x_1,\ldots,x_n)\in X_1\times\ldots\times X_n$ in directions $\eta_i\in X_i$, $1\leq i \leq k$, is
\begin{equation*}
\delta^k_{1:k} F \left( x; \eta_1,\ldots,\eta_k \right) = F(\eta_1,\ldots,\eta_k,x_{k+1},\ldots,x_n),
\end{equation*}
where the differential $\delta^k_{1:k}F \left( x; \eta_1,\ldots,\eta_k \right)$ refers to $\delta_k ( \ldots \delta_2(\delta_1 F ( x; \eta_1 );\eta_2) \ldots ;\eta_k)$.

The $k^{th}$-order total chain differentials of $F$ at point $x=(x_1,\ldots,x_n)\in X_1\times\ldots\times X_n$ in directions $\underline{\eta}_i\in X_1\times\ldots\times X_n$, $1\leq i \leq k$, is
\begin{equation*}
\delta^k F \left( x; \underline{\eta}_1,\ldots,\underline{\eta}_k \right) = \sideset{}{^{\neq}}\sum_{1\leq j_1,\ldots,j_k \leq n} F\left(y_{j_{1:k}}^1,\ldots,y_{j_{1:k}}^n\right),
\end{equation*}
where if there exists $r\in\mathbb{N}$, $1\leq r \leq k$, such that $i = j_r$ then $y_{j_{1:k}}^i = \underline{\eta}^i_r$ and $y_{j_{1:k}}^i = x_i$ otherwise. The sum $\sum^{\neq}$ is over different indices $j_1\neq\ldots\neq j_k$.
\end{mylemma}

\begin{IEEEproof}
The base case $k=0$ is obvious. Assuming the relation is true at order $k-1$, $k^{th}$-order chain differentials are
\begin{align*}
\delta^k_{1:k} F \left( x; \eta_1,\ldots,\eta_k \right) = \delta_k \left( F(\eta_1,\ldots,\eta_{k-1},x_k,\ldots,x_n) ;\eta_k\right)
\end{align*}
and
\begin{align*}
\delta^k F (x; \underline{\eta}_1,\ldots,\underline{\eta}_k) & = \sideset{}{^{\neq}}\sum_{1\leq j_1,\ldots,j_{k-1} \leq n} \, \delta \left(F\left(y_{j_{1:k-1}}^1,\ldots,y_{j_{1:k-1}}^n\right);\underline{\eta}_k\right) \\
& = \sideset{}{^{\neq}}\sum_{1\leq j_1,\ldots,j_{k-1} \leq n} \, \sum_{\substack{j_k=1\\j_k\neq j_{1:k-1}}}^n F\left(y_{j_{1:k}}^1,\ldots,y_{j_{1:k}}^n\right).
\end{align*}
Both equations are according to the claimed result at order $k$. The lemma is then proved by induction.
\end{IEEEproof}

The next section describes stochastic population processes, and uses the G\^{a}teaux differential to determine factorial moment and probability measures from probability generating functionals.

\section{Stochastic population processes}
\label{sec:stoPopProc}

Moyal \cite{Moyal62} developed the general theory of stochastic population processes, as a generalisation of stochastic processes for modelling collections of individual state variables.
These were based on generating functionals~\cite{Volterra}, which are generalised versions of generating functions~\cite{deMoivre}.
Moyal's general theory of stochastic population processes provides the mathematical foundation for point process theory~\cite{DaleyVere-Jones,RCoxIsham} and branching processes~\cite{harris}.
The probability and factorial moment measures are found with G\^ateaux differentials.

We extend the probability generating functional to permit non-symmetric joint measures, which is achieved by relinquishing the commutativity of the product in the formal power series that describes the generating functional.
The consequence of this is that it allows to work with point processes with non-symmetric measures.
Non-symmetric point processes were originally proposed by Moyal \cite{Moyal62}, yet the commutativity of the product in the polynomial did not permit such processes to be defined with the usual generating functional of Volterra \cite{Volterra}.
The purpose of this generalisation is that it allows us to unify the problems of branching processes, renormalization, and more general interacting processes, within the same paradigm.

\begin{mydef}[Stochastic population processes]
A stochastic population process is a measurable mapping $\varphi:(\Omega,\mathcal{F},P)\to(\mathcal{X},\mathbf{B}_{\mathcal{X}})$ for a given probability space $(\Omega,\mathcal{F},P)$ and measurable space $(\mathcal{X},\mathbf{B}_{\mathcal{X}})$.

The process state space $\mathcal{X}$ is the space of all the sets of points in $\mathbf{X}=\mathbb{R}^d$ and $\mathbf{B}_{\mathcal{X}}$ is the Borel $\sigma$-algebra on $\mathcal{X}$.
The probability measure $P$ in the space $\mathcal{P}(\mathcal{X})$ of all probability measures on $(\mathcal{X},\mathbf{B}_{\mathcal{X}})$ depicts the distribution of $\varphi$.
\end{mydef}

\subsection{Probability generating functional}

The probability generating functional (p.g.fl.) provides the mathematical foundation for point process theory~\cite{Moyal62}, which is defined as follows

\begin{mydef}[Probability generating functional]
Let $w^{s}$ be a symmetric function in $\mathcal{B}_b(\mathcal{X})$, the set of bounded measurable function on $(\mathcal{X},\mathbf{B}_{\mathcal{X}})$, and $w^{s}(x_{1:n}) = w(x_1)\ldots w(x_n)$. The probability generating functional of a process $\varphi$ with probability measure $P\in\mathcal{P}(\mathcal{X})$ is defined with $G(w) = \mathbb{E}[ w^{s}(\varphi)]$ or, written differently
\begin{align}
\label{eq:defPgfl}
G(w) & = \int w^{s}(\varphi) dP \\
& = \sum_{n\geq 0} \int w(x_1)\ldots w(x_n) P(dx_1\times\ldots\times dx_n).
\end{align}
\end{mydef}

Since the test function $w^{s}$ is symmetric, two countably equivalent distributions have the same p.g.fl.~\cite{Moyal62}. Accounting for non-symmetric probability measures is then not possible.
To account for processes with non-symmetric probability measures, it is necessary to generalize (\ref{eq:defPgfl}) for non-symmetric functions.

\begin{mydef}[Non-symmetric test function]
Let $\mathbf{X}$ be the individual state space and let us consider a functional $F:\mathcal{B}_b(\mathbf{X})\to\mathcal{B}_b(\mathcal{X})$ and a function $w\in\mathcal{B}_b(\mathbf{X})$. Then the function $F(w)$ is defined, for all $x_i \in \mathbf{X}$, for all $n\in\mathbb{N}$, by
\begin{align}
F(w)(x_1,\ldots,x_n) & = F(w,\ldots,w)(x_1,\ldots,x_n)\\
& = w(x_1)\ldots w(x_n).
\end{align}
\end{mydef}

\begin{mycorollary}[Partial and total chain differential of $F(w)$]
\label{corol:wDiff}
The $k^{th}$-order partial and total chain differentials of $F(w)(x_{1:n})$, $n\in\mathbb{N}$, in directions $\eta_i\in\mathcal{B}_b(\mathbf{X})$ and $\underline{\eta}_i\in\mathcal{B}_b(\mathbf{X})^n$, $1\leq i \leq k$ are respectively
\begin{align*}
\delta^k_{1:k} \left( F(w)(x_{1:n}); \eta_1,\ldots,\eta_k \right) & = \prod_{i=1}^k\eta_i(x_i)\prod_{i=k+1}^n w(x_i),\\
\delta^k \left( F(w)(x_{1:n}); \underline{\eta}_1,\ldots,\underline{\eta}_k \right) & =
\sideset{}{^{\neq}}\sum_{1\leq j_1,\ldots,j_k \leq n} \, \prod_{i=1}^n \mu^i_{j_{1:k}}(x_i)
\end{align*}
where if there exists $r\in\mathbb{N}$, $1\leq r \leq k$, such that $i = j_r$, then $\mu^i_{j_{1:k}} = \underline{\eta}^i_r$ and $\mu^i_{j_{1:k}} = w$ otherwise.
\end{mycorollary}

For $\mathbb{E}[F(w)(\varphi)]$ to be a convenient representation of the process, one has to be able to recover probability measures and factorial moment measures from it. For some applications, the symmetrized form $J$ of the probability measure $P$, named Janossy measure \cite{JanossyCascade,DaleyVere-Jones}, also has useful properties and is defined by
\begin{equation}
J(dx_1\times\ldots\times dx_n) = \sum_{\sigma(1:n)} P(dx_{\sigma_1}\times\ldots\times dx_{\sigma_n})
\end{equation}
where $\sigma(1:n)$ is the set of all permutations of the set of indices \{1:n\}. The objective of the following corollary is to recover these three different ways of characterizing the stochastic process $\varphi$ from the expectation $\mathbb{E}[F(w)(\varphi)]$.

\begin{mycorollary}[From Corollary \ref{corol:wDiff} and Lemma \ref{lemma:expectDiffInvers}]
\label{corol:RecoverPPdescriptor}
Let $\mathcal{M}(\mathcal{X})$ be the set of measures on $(\mathcal{X},\mathbf{B}_{\mathcal{X}})$. The way of recovering the probability measure $P\in\mathcal{P}\left(\mathcal{X}\right)$, the Janossy measure $J\in\mathcal{M}\left(\mathcal{X}\right)$ and the factorial moment measure $M\in\mathcal{M}(\mathcal{X})$ from the expectation $\mathbb{E}\left[ F(w)(\varphi)\right]$ is given by
\begin{align*}
P(dx_{1:k}) & = \left. \delta^k_{1:k} \left( \mathbb{E}\left[ F(w)(\varphi)\right]; \mathbf{1}_{dx_1},\ldots,\mathbf{1}_{dx_k} \right) \right|_{w = 0},\\ 
J(dx_{1:k}) & = \left. \delta^k \left( \mathbb{E}\left[ F(w)(\varphi) \right]; \mathbf{1}_{dx_1},\ldots,\mathbf{1}_{dx_k} \right) \right|_{w = 0},\\
M(dx_{1:k}) & = \left. \delta^k \left( \mathbb{E}\left[ F(w)(\varphi) \right]; \mathbf{1}_{dx_1},\ldots,\mathbf{1}_{dx_k} \right) \right|_{w = 1},
\end{align*}
where $\mathbf{1}_{dx}$ is the indicator function of $dx$.
\end{mycorollary}

Intuitively speaking, the $k^{th}$-order partial differentiation at $w=0$ of $\mathbb{E}[ F(w)(\varphi)]$ in the directions $\{\mathbf{1}_{dx_i}\}_{i=1}^k$ imposes that there are $k$ points in the process $\varphi$ and the $i^{th}$ point is in the infinitesimal region $dx_i$. The $k^{th}$-order total differentiation at $w=0$ imposes that there are $k$ points with exactly one point in each $dx_i$. The $k^{th}$-order total differentiation at $w=1$ only imposes that any $k$ points in the process are in the regions $\{dx_i\}_{i=1}^k$. It is interesting to note that the total differentiation of $F(w)(x_{1:n})$ is equal to the differentiation of $w(x_1)\ldots w(x_n)$. Therefore, the Janossy measure $J$ and the factorial moment $M$ can be recovered equally from the associated p.g.fl. $G(w)$, but not the probability $P$.

\begin{mydef}[The counting measure $C_k$] \label{def:countingMeasure} The total chain differential $\delta^k \left. \left(F(w);\mathbf{1}_{dx_1},\ldots,\mathbf{1}_{dx_k} \right) \right|_{w = 1}$ is a counting measure denoted $C_k(dx_{1:k}|\cdot)$.
\end{mydef}

The counting measure $C_k$ is equivalent to the counting measure $N_{(k)}$ defined by Moyal \cite{Moyal62}, p.11, eq. (3.11).

\begin{myexample}[Independent stochastic process]
The probability measure $P\in\mathcal{P}(\mathcal{X})$ of an independent stochastic process $\varphi$ is defined by its projection on $\mathbf{X}^n$:
\begin{equation*}
p(n) F(P^{(n)}_1,\ldots,P^{(n)}_n)
\end{equation*}
so that, for all $n\in\mathbb{N}$, $(x_1,\ldots,x_n) \in \mathbf{X}^n$,
\begin{equation*}
P(dx_1\times\ldots\times dx_n) = p(n) \prod_{i=1}^n P^{(n)}_i(dx_i).
\end{equation*}
Such a process can be easily described by the quantity $\mathbb{E}[F(w)(\varphi)]$ whereas the usual probability generating functionals could only generate the symmetrised form
\begin{equation*}
\dfrac{1}{n!}\sum_{\sigma} P^{(n)}_i(dx_{\sigma_i}),
\end{equation*}
for any given $n\in\mathbb{N}$.

It is commonly assumed in multi-target tracking applications that each target is independent of each other and has its own target distribution. This can be represented by a multi-Bernoulli process \cite{MahlerBook} with parameter set
\begin{equation*}
\left\{ \{q_{n,i},P^{(n)}_i\}_{i=1}^n \right\}_{n\geq 0}
\end{equation*}
and its probability measure reads
\begin{equation*}
P(dx_1\times\ldots\times dx_n) = \prod_{i=1}^n q_{n,i} P^{(n)}_i(dx_i).
\end{equation*}

Multi-Bernoulli processes can be seen as independent stochastic processes where $p(n) = \prod_{i=1}^n q_{n,i}$.
\end{myexample}

Since we shall be dealing with bivariate processes to define hierarchy and determine the Chapman-Kolmogorov equation and Bayes' rule, we define the corresponding bivariate generalised stochastic population process as follows.

\begin{mydef}[Joint probability generating functional]
Let us consider the joint process $\varphi = (\varphi_x, \varphi_y)$ on the measurable space $\left(\mathcal{X}\times\mathcal{Y},\mathbf{B}_{\mathcal{X}}\otimes\mathbf{B}_{\mathcal{Y}}\right)$. Let us then define a functional $F:(w_x,w_y)\in\mathcal{B}_b(\mathbf{X}) \times \mathcal{B}_b(\mathbf{Y}) \mapsto F_x(w_x)F_y(w_y)\in\mathcal{B}_b(\mathcal{X}) \times \mathcal{B}_b(\mathcal{Y})$, and the test function $w=(w_x,w_y)$ where $w_x \in \mathcal{B}_b(\mathbf{X})$ and $w_y \in \mathcal{B}_b(\mathbf{Y})$. The expectation of $F(w)(\varphi)$ with probability measure $P_{x,y} \in \mathcal{P}(\mathcal{X} \times \mathcal{Y})$ is
\begin{align*}
\mathbb{E}\left[ F(w)(\varphi) \right] & = \int F(w_x,w_y)(\varphi_x,\varphi_y) dP\\
& = \sum_{n,m\geq 0} \int F_x(w_x)(x_{1:n})F_y(w_y)(y_{1:m}) P_{x,y}(dx_{1:n}\times dy_{1:m}).
\end{align*}
\end{mydef}

In the next section, the tools developed here for describing stochastic population processes are used along with Fa\`a di Bruno's formula for chain differentials (Theorem \ref{thm:chainRule}) to find the probability measures and the factorial moment measures of hierarchical processes.

\section{Interacting stochastic population processes}
\label{sec:hierInterPopProc}

In this section we introduce a new class of stochastic population processes for modelling interactions and hierarchies.
Firstly, we discuss hierarchical population processes, that are constructed by composing expectations.
Secondly, we introduce interacting population processes, that are able to model dependencies between groups of parent and daughter processes.
Then, this general population process enable us to unify branching processes and renormalization within the same paradigm in the next two sections.

To define interactions, the notion of Markov kernel \cite{pdm} is needed.

\begin{mydef}[Markov kernel]
A Markov kernel $K$ from a measurable space $(\mathcal{X},\mathbf{B}_{\mathcal{X}})$ into another measurable space $(\mathcal{Y},\mathbf{B}_{\mathcal{Y}})$ is an integral kernel such that $K(\cdot|x)$ is in $\mathcal{P}(\mathcal{Y})$ for all $x\in\mathcal{X}$ and the mapping $x \xmapsto{} K(A_y|x)$ is a $\mathcal{X}$-measurable function for all $A_y\in\mathbf{B}_{\mathcal{Y}}$.

The set of Markov kernels from $(\mathcal{X},\mathbf{B}_{\mathcal{X}})$ to $(\mathcal{Y},\mathbf{B}_{\mathcal{Y}})$ is denoted $\mathcal{K}(\mathbf{B}_{\mathcal{Y}}\times\mathcal{X})$.
\end{mydef}

\subsection{Hierarchical population processes}

Hierarchical population processes form a general class of process composed of two processes named parent and daughter processes. The daughter process is conditioned on the parent process. Unlike cluster processes, daughters of hierarchical processes are conditioned on the whole parent process. Cluster processes, branching processes, data clustering or multi target tracking are instances of hierarchical processes. Also, the Chapman Kolmogorov equation and Bayes' rule can be easily formulated in terms of hierarchical population processes between the states at two different times and between measurements and state respectively.

\begin{mytheorem}[Hierarchical population processes]
\label{thm:hierarchicalProcess}
Let $\varphi_{p,d}=(\varphi_p,\varphi_d)$ be a joint process on $\left(\mathcal{X}\times\mathcal{Y},\mathbf{B}_{\mathcal{X}}\otimes\mathbf{B}_{\mathcal{Y}}\right)$ composed of parent and daughter processes $\varphi_p$ and $\varphi_d$ on $(\mathcal{X},\mathbf{B}_{\mathcal{X}})$ and $(\mathcal{Y},\mathbf{B}_{\mathcal{Y}})$ respectively.
The expectation w.r.t. the joint parent/daughter process is
\begin{equation*}
\mathbb{E}\left[F_{d,p}(w_d,w_p)(\varphi_{p,d})\right] = \mathbb{E}\left[F_p(w_p)(\varphi_p) \mathbb{E}\left[\left.F_d(w_d)(\varphi_d)\right|\varphi_p\right]\right].
\end{equation*}
where $\mathbb{E}\left[\,\cdot\,|\varphi_p\right]$ is the conditional expectation of the daughter process given the parent process.
\end{mytheorem}

\begin{IEEEproof}
Let us consider the decomposition of the joint probability measure $P_{p,d}\in\mathcal{P}(\mathcal{X}\times\mathcal{Y})$ into the Markov kernel $P_{d|p}\in\mathcal{K}(\mathbf{B}_{\mathcal{Y}}\times\mathcal{X})$ and the probability measure $P_p$ as $P_{p,d}\left( dx_{1:k},dy_{1:l} \right) = P_{d|p}\left( dy_{1:l}|x_{1:k}\right)P_p\left( dx_{1:k} \right)$. Considering the expectation of $F_{p,d}(w_p,w_d)(\varphi_{p,d})$ where the probability measure of $\varphi_{p,d}$ is $P_{p,d}$ gives the desired result.
\end{IEEEproof}

It is often necessary to write the conditional expectation of $\varphi_d$ given $\varphi_p$ in a different form to take into account the possible assumptions made on the hierarchical process to define other more specific processes. The main assumptions needed are about independence since hierarchical processes have a high degree of internal correlations. Interacting processes (Section \ref{ssec:hierInterPopProc}) is a powerful tool when it comes to designing the process $\varphi_{p,d}$ given some independence properties.

The objective of the next section is to provide a useful framework for designing the conditional expectation of $\varphi_d$ given $\varphi_p$.

\subsection{Interacting population processes}
\label{ssec:hierInterPopProc}

In this section we develop models for representing hierarchical systems of objects inspired from connected correlation functions~\cite{ursell}.
The models and examples developed will be used for physical models in the following sections.

Let us consider a process $\psi$ on $(\mathcal{Z},\mathbf{B}_{\mathcal{Z}})$ and two joint processes $\phi_{p,d}$ and $\varphi_{p,d}$ on $\left(\mathcal{X}\times\mathcal{Y},\mathbf{B}_{\mathcal{X}}\otimes\mathbf{B}_{\mathcal{Y}}\right)$. The process $\varphi_{p,d}$ is conditioned on one point $z$ in the parameter space $\mathbf{Z}$.
We consider the process $\phi_{p,d}$ to be resulting from the superposition of a random number of instances of the process $\varphi_{p,d}$ according to another process $\psi$ in the following way:
\begin{equation}
\label{eq:jointForInteractingKernels}
\mathbb{E}[F_{p,d}(w_{p,d})(\phi_{p,d})] = \mathbb{E}_{\psi}[ F_{\psi}(\mathbb{E}[F_{p,d}(w_{p,d})(\varphi_{p,d})|\cdot])(\psi)].
\end{equation}

The main objective is to deduce the relation between the superposed process $\phi_d$ and the simple process $\varphi_d$ composing the superposition.

\begin{mytheorem}[Interacting population processes]
\label{thm:hierInteractingPopProc}
The process $\phi_d$ generated from (\ref{eq:jointForInteractingKernels}) is described through the following conditional expectation
\begin{equation*}
\mathbb{E}[F_d(w_d)(\phi_d) |\phi_p = x_{1:k}] \propto \sum_{\pi \in \Pi(x_{1:k})} \mathbb{E}_{\psi}\left[ Q(F_d,w_d,\pi,\psi) \right],
\end{equation*}
where
\begin{equation*}
Q(F,w,\pi,z_{1:n}) = \prod_{i=1}^{\dot{\pi}} p_p(\dot{\pi}_i)\mathbb{E} \left[ F(w)(\varphi_d) | \pi_i,z_i \right]\prod_{i=\dot{\pi}+1}^n p_p(0)\mathbb{E}\left[ F(w)(\varphi_d) | \emptyset, z_i \right].
\end{equation*}

The normalization is found by expanding $F_d(w_d)$ as a product and then considering $w = 1$.
Note that the probability mass function (p.m.f.) $p_p$ could be considered as dependent on the points $z_i$ of the process $\psi$.
\end{mytheorem}

\begin{IEEEproof}
Let us first consider the expectation
\begin{equation*}
\mathbb{E}\left[ F_{p,d}(w_{p,d})(\varphi_{p,d}) \right] = \mathbb{E} \left[ F_p(w_p)(\varphi_p) \mathbb{E}\left[ F_d(w_d)(\varphi_d) | \varphi_p \right] \right].
\end{equation*}
Since we are interested in the conditional expectation of $\varphi_d$ given $\varphi_p$, we assume the spatial distribution of $\varphi_p$ is known, but not the associated p.m.f. $p_p$, so the probability measure is $P_p(dx_{1:n}) = p_p(n)\mathbf{1}_{dx_1}(x_1),\ldots,\mathbf{1}_{dx_n}(x_n)$, for all $n\in\mathbb{N}$, where $x_i \in \varphi_p$. The way to recovering the conditional expectation of $F_d(w_d)(\varphi_d)$ given $\varphi_p$ from the expectation of $F_{p,d}(w_{p,d})(\varphi_{p,d})$ is provided by the following chain differential:
\begin{equation*}
\mathbb{E}[ F_d(w_d)(\varphi_d) | \varphi_p = x_{1:k}] = \dfrac{1}{p_p(k)}\left.\delta^k_{1:k} \mathbb{E}\left[ F_{p,d}(w_{p,d})(\varphi_{p,d}); 1,\ldots,1 \right]\right|_{w_p = 0}.
\end{equation*}

Let us now consider the expectation of $\phi_{p,d}$ defined in (\ref{eq:jointForInteractingKernels}). The conditional expectation of $F_d(w_d)(\phi_d)$ given $\phi_p$ can then be recovered  using Fa\`a di Bruno's formula for chain differentials and then considering the probability measures $P_p$ and $\bar{P}_p$ to be
\begin{align*}
P_p(dx_{1:n}) & = p_p(n)\mathbf{1}_{dx_1}(x_1),\ldots,\mathbf{1}_{dx_n}(x_n)\\
\bar{P}_p(dx_{1:n}) & = \bar{p}_p(n)\mathbf{1}_{dx_1}(x_1),\ldots,\mathbf{1}_{dx_n}(x_n).
\end{align*}

The p.m.f. $p_p$ is considered as a model parameter while $\bar{p}_p$ is an unknown constant to be determined be considering $w_d = 1$ in the final result.

\end{IEEEproof}

Henceforth, the indices of $F$ and $w$ are going to be omitted when there is no ambiguity.

From here, designing a transition just consists in choosing a specific p.m.f. $p_p$ and a distribution for $\psi$. To find the probability measure related to interacting stochastic population processes, it is possible to consider either the partial chain differential of $\mathbb{E}[F(w)(\varphi_d) | x_{1:k}]$ w.r.t. the daughter process or the partial chain differential of (\ref{eq:jointForInteractingKernels}) w.r.t. both parent and daughter processes. The former gives a detailed but involved expression while the latter gives a more general and simpler overview of the process and is therefore preferred.

\begin{mytheorem}
\label{thm:ProbaAndFactOfHierIntPopProc}
The Markov kernel $\bar{P}_{d|p}\in\mathcal{K}(\mathbf{B}_{\mathcal{Y}}\times\mathcal{X})$ and the associated factorial moment measure $\bar{M}_{d|p}$ depicting the probability of $\phi_d$, given that the realisation of the parent process $\varphi_p$ is $\{x_{1:k}\}$, are
\begin{align*}
\bar{P}_{d|p}\left( dy_{1:l} | x_{1:k} \right) &\propto \sum_{\pi \in \Pi(y_{1:l}\cup x_{1:k})} \mathbb{E}_{\psi} \left[Q_P(\pi,\psi)\right]\\
\bar{M}_{d|p}\left( dy_{1:l} | x_{1:k} \right) &\propto \sum_{\pi \in \Pi(y_{1:l}\cup x_{1:k})} \mathbb{E}_{\psi} \left[Q_M(\pi,\psi)\right],
\end{align*}
with
\begin{align*}
Q_P(\pi,z_{1:n}) & = \prod_{i=1}^{\dot{\pi}} p_p(\dot{\pi}_i)P_{d|p}\left(d\pi_{i,y} | \pi_{i,x},z_i\right)\prod_{i=\dot{\pi}+1}^n p_p(0)P_{d|p}\left(\emptyset|\emptyset,z_i\right),\\
Q_M(\pi,z_{1:n}) & = \sideset{}{^{\neq}}\sum_{1 \leq j_1,\ldots,j_{\dot{\pi}} \leq n} p_p(0)^{n-\dot{\pi}} \prod_{i=1}^{\dot{\pi}} p_p(\dot{\pi}_i) M_{d|p}\left(d\pi_{j_i,y}|\pi_{j_i,x},z_i\right),
\end{align*}
where $\pi_{i,x}$ stands for $\pi_i \cap \{x_{1:k}\}$. The normalizing constant is the same as in Theorem \ref{thm:hierInteractingPopProc}.
\end{mytheorem}

\begin{IEEEproof}
The Markov kernel $\bar{P}_{d|p}$ and the factorial moment measure $\bar{M}_{d|p}$ are recovered from (\ref{eq:jointForInteractingKernels})
respectively by partial and total differentiation (see Corollary \ref{corol:RecoverPPdescriptor}):
\begin{align*}
\bar{P}_{d|p}(dy_{1:l}|x_{1:k}) & = \delta^k_{1:k}(\delta^l_{1:l} ( \mathbb{E}\left[ F(w)(\phi_{p,d}) \right];\mathbf{1}_{dy_1},\ldots,\mathbf{1}_{dy_l}) ;1,\ldots,1 )|_{w=(0,0)},\\
\bar{M}_{d|p}(dy_{1:l}|x_{1:k}) & = \delta^{l+k} ( \mathbb{E}\left[ F(w)(\phi_{p,d}) \right];\mathbf{1}_{dy_1},\ldots,\mathbf{1}_{dy_l},\underline{1}',\ldots,\underline{1}')|_{w=(0,1)},
\end{align*}
where $\underline{1}' = [1,0,\ldots,0]^T$. The quantity $\underline{1}'$ enables to write partial chain differentials as total differentials and therefore to use a $(l+k)^{th}$-order total chain differential to find the factorial moment measure $\bar{M}_{d|p}$. The calculation of the moment involves terms equal to $\delta^0 (\mathbb{E}[ F(w)(\phi_{p,d})])|_{w=1} = 1$.
\end{IEEEproof}

\begin{myexample}[Multiplicative population processes]
Moyal \cite{MoyalAP,Moyal62,MoyalRS} introduced the concept of multiplicative population processes to describe a Markov process where individuals at some time step $s$ are the ``ancestors'' of mutually independent populations at time $t\ge s$.
This can be characterised by the superposition of point processes, each conditioned on a single state.
The process is described with the following expression:
\begin{equation*}
\mathbb{E}\left[ F(w)(\phi_d) | \varphi_p = x_{1:k} \right] = \prod_{i=1}^k \mathbb{E}\left[ F(w)(\varphi_d) | x_i \right],
\end{equation*}
which is obtained when considering in Theorem \ref{thm:hierInteractingPopProc} that
\begin{itemize}
\item $\varphi_d$ is independent of $\psi$,
\item the p.m.f. $p_{\psi}(n) = \int P_{\psi}(dx_{1:n})$ is equal to one if $n=k$ and zero otherwise, and
\item the p.m.f. $p_p(n)$ is equal to one at $n=1$ and zero otherwise.
\end{itemize}
These assumptions make each process in the superposition only dependent on one individual of the parent process $\varphi_p$, see Figure \ref{fig:clusterProcess}.

\begin{figure}[htp]
\centering\includegraphics[width=300px]{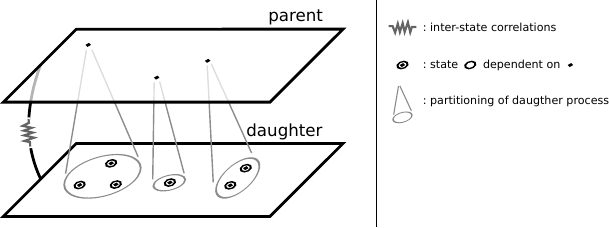}
\caption{Transition of a multiplicative population process.}
\label{fig:clusterProcess}
\end{figure}
\end{myexample}

\subsection{Superpositional processes and Khinchin processes}
\label{ssec:SuperProcAndKhinchinProc}

In this section we describe two related classes of processes. The first one, named superpositional process, is a superposition of stochastic population processes.
The second one, named Khinchin process, is a stochastic population process with probability measure based on Khinchin measures.
First, each of these two classes of processes are described. Then, their relation and a application example are given.

\begin{mydef}[Superpositional processes]
Let $\phi$ and $\varphi$ be processes on $(\mathcal{X},\mathbf{B}_{\mathcal{X}})$ and $\psi$ be another process on $(\mathcal{Z},\mathbf{B}_{\mathcal{Z}})$. The process $\varphi$ is conditioned on one point $z$ in the parameter space $\mathbf{Z}$.
Let us consider $\phi$ to be generated by the superposition of a random number of instances of $\varphi$ according to $\psi$ such that
\begin{equation}
\label{eq:superposProc}
\mathbb{E}\left[ F(w)(\phi) \right] = \mathbb{E}_{\psi}\left[ F_{\psi}(\mathbb{E}\left[ F(w)(\varphi) | \cdot \right])(\psi) \right],
\end{equation}
then $\phi$ is said to be a superpositional process.
\end{mydef}

\begin{myexample}[Poisson superpositional processes]
Let $\varphi$ be independent of $\psi$ and $p_{\psi}(k) = k!^{-1} \lambda^k \exp(-\lambda)$. The expectation (\ref{eq:superposProc}) can be rewritten
\begin{equation*}
\mathbb{E}\left[ F(w)(\phi) \right] = \sum_{k\geq 0} \dfrac{1}{k!} \lambda^k e^{-\lambda} f_{\psi}(\mathbb{E}\left[ F(w)(\varphi)\right](k),
\end{equation*}
where $f_{\psi}(w)(k) = w^k$.
This form is not particularly convenient and does not make the Poisson assumption useful. But if we expand $f_{\psi}$ and $F$ to simple products, this expression becomes
\begin{equation*}
\mathbb{E}\left[ w^s(\phi) \right] = \exp\left( \lambda\mathbb{E}\left[ w^s(\varphi) - 1\right]\right).
\end{equation*}

However, a part of the structure has been lost by doing so, and the consequence is that only symmetrized quantities as Janossy or factorial moment measures can be recovered.

\end{myexample}

What have been demonstrated for Poisson superpositional processes is true for every process based on exponential distributions.
To recover the Janossy and factorial moment measures describing this kind of process, the differentiation is w.r.t. $w$ and not $F$ as earlier.

\begin{mydef}[Khinchin processes]
Let $\{K_n\}_{n\geq1}$ be a sequence of Khinchin measures such that for all $n\in\mathbb{N}\setminus \{0\}$, $K_n \in \mathcal{M}(\mathbf{X}^n)$ and let $\phi$ be a process on $(\mathcal{X},\mathbf{B}_\mathcal{X})$ defined through the expectation $\mathbb{E}\left[ w^{s}(\phi) \right]$ by
\begin{equation*}
\mathbb{E}\left[ w^{s}(\phi)\right] = \exp\left( -K_0 + \sum_{n\geq1} \int w^{s}(x_{1:n}) K_n(dx_{1:n}) \right),
\end{equation*}
where, when denoting $K_n(\mathbf{X}^n)$ the integral $\int K_n(dx_{1:n})$, $K_0$ is equal to $\sum_{n\geq1} K_n(\mathbf{X}^n)$.
The process $\phi$ is named Khinchin process.
\end{mydef}

\begin{myexample}[Poisson Khinchin process]
Let $K_n$ be a Poisson intensity measure $\lambda_n \in \mathcal{M}(\mathbf{X}^n)$, the Khinchin process $\phi$ becomes
\begin{equation*}
\mathbb{E}\left[ w^{s}(\phi) \right] = \exp\left(\sum_{n\geq1} \int \left(w^{s}(x_{1:n}) - 1 \right) \lambda_n(dx_{1:n}) \right),
\end{equation*}
which is the generalisation of Poisson and Gauss-Poisson processes justifying the name Poisson Khinchin process.
\end{myexample}

\begin{myexample}[Khinchin processes as a superpositional processes]
Let $K_n$ be the projection of $P\in\mathcal{P}(\mathcal{X})$ on the space $\mathbf{X}^n$, the Khinchin process $\phi$ becomes
\begin{equation*}
\mathbb{E}\left[ w^{s}(\phi) \right] = \exp\left( \sum_{n\geq0} \int w^{s}(x_{1:n}) P_n(dx_{1:n}) - 1 \right),
\end{equation*}
since $\sum_{n\geq1}P_n(\mathbf{X}^n) = 1 - P_0$.

We conclude that a Khinchin process is a Poisson superpositional processes (with parameter $\lambda = 1$ and $F(w)$ developed as a product) when its measures $K_n$ are probability measures.
\end{myexample}

\begin{myexample}[Jointly Khinchin processes]
\label{ex:JointKhinchin}
Let us consider a jointly Khinchin processes $\varphi$ and $\phi$ on $(\mathcal{X},\mathbf{B}_{\mathcal{X}})$ and $(\mathcal{Y},\mathbf{B}_{\mathcal{Y}})$ respectively, with Khinchin measures $K_{n,m} \in \mathcal{M}(\mathbf{X}^n\times\mathbf{Y}^m)$ such that
\begin{equation*}
K_{n,m}(dx_{1:n}\times dy_{1:m}) = P(dx_{1:n}|y_{1:m})\lambda_m(dy_{1:m})
\end{equation*}
where $P\in\mathcal{K}(\mathbf{B}_{\mathcal{X}}\times\mathcal{Y})$ is a Markov kernel and $\lambda_m\in\mathcal{M}(\mathbf{Y}^m)$ is Poisson measure. The expectation $\mathbb{E}[w^{s}(\varphi,\phi)]$ can be written
\begin{equation}
\label{eq:JointKhinchin}
\mathbb{E}[w^{s}(\varphi,\phi)] = \exp\left( H(w^s) - H(1)\right),
\end{equation}
where, assuming $\lambda_0 = 1$,
\begin{equation*}
H(w^s) = \sum_{\substack{n,k\geq 0 \\ (n,k)\neq(0,0)}} \int w^s(x_{1:n},y_{1:k}) P(dx_{1:n}|y_{1:k})\lambda_k(dy_{1:k}).
\end{equation*}

Differentiating (\ref{eq:JointKhinchin}) w.r.t. $w_x$ and $w_y$ in directions $\{\mathbf{1}_{dx_1},\ldots,\mathbf{1}_{dx_n}\}$ and $\{\mathbf{1}_{dy_1},\ldots,\mathbf{1}_{dy_m}\}$ and at points $w_x=0$ and $w_y=1$, we find
\begin{equation*}
J(dx_{1:n}|y_{1:m})M(dy_{1:m}) \propto
\sum_{\pi \in \Pi(x\cup y)} \prod_{\omega\in\pi} \sum_{k\geq 0} \int \sideset{}{^{\neq}}\sum_{1\leq j_1,\ldots,j_{\dot{\omega}_y} \leq k} J(d\omega_x|\hat{y}^{j_{1:\dot{\omega}_y}}_{1:k})\lambda_k(d\hat{y}^{j_{1:\dot{\omega}_y}}_{1:k}).
\end{equation*}
where for all $r\in\mathbb{N}$, $1\leq r \leq k$, $\displaystyle\hat{y}^{j_{1:\dot{\omega}_y}}_r$ is defined to be
\begin{align}
\left\{
\begin{array}{ll}
\omega_{y,i} & \text{ if } \exists i\in\mathbb{N}, 1\leq i\leq\dot{\omega}_y, r = j_i \\
\bar{y}_r & \text{ otherwise}.
\end{array}
\right.
\end{align}

This result, along with the Bayes' rule, provides a generalisation of the Probability Hypothesis Density filter (first order, see \cite{Mahlermaths}) and of filters with a Gauss-Poisson prior (second order, see Bochner \cite{Bochner}, p. 145) to any order, as demonstrated in the sequel.
\end{myexample}

\section{The Chapman-Kolmogorov equation}
\label{sec:ChapmanKolmogorov}

\begin{figure}[htp]
\centering
\includegraphics[width=300px]{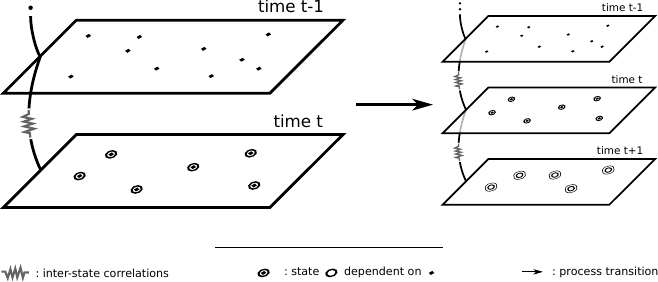}
\caption{Transition of a point process.}
\label{fig:ChapmanKolmogorov}
\end{figure}

The Chapman-Kolmogorov equation was first proposed by Einstein \cite{Einstein} to describe the motion of interacting molecules in his theory of Brownian movement.
This concept was later investigated by Chapman \cite{chapman} for displacements of grains suspended in non-uniform fluid, and Kolmogorov \cite{kolmogorov} in probability theory.
For a Markovian stochastic process, the Chapman-Kolmogorov equation depicts the joint probability distributions after applying a transition probability on the process (see Figure \ref{fig:ChapmanKolmogorov}). Times indices are traditionally used to describe the different states involved in the Chapman-Kolmogorov equation though this equation describes any kind of transition between two state spaces provided that
a Markov transition is defined between them.
In this section we describe the Chapman-Kolmogorov equation for systems of multiple objects with interactions and hierarchies.

The following corollary provides the formal description of the Chapman-Kolmogorov equation that we are interested in.
\begin{mycorollary}[From Theorem \ref{thm:hierarchicalProcess}]
\label{corol:ChapmanKolmogorov}
Let $\varphi_{t+1}$ and $\varphi_t$ be processes with respective measures $P_{t+1}\in\mathcal{P}\left(\mathcal{X}\right)$, $P_t \in \mathcal{P}\left(\mathcal{Y}\right)$ and let $P_{t+1|t}\in\mathcal{K}(\mathbf{B}_{\mathcal{X}}\times\mathcal{Y})$ be a Markov kernel. The Chapman-Kolmogorov equation for multi-object processes is
\begin{equation*}
\mathbb{E}\left[F(w)(\varphi_{t+1})\right] = \mathbb{E}\left[\mathbb{E}\left[F(w)(\varphi_{t+1})|\varphi_t\right]\right]
\end{equation*}
where $F(w)\in\mathcal{B}_b(\mathcal{X})$ is a test function.
\end{mycorollary}

The Chapman-Kolmogorov equation can be interpreted as a hierarchical process in time, where the parent process represents the previous state, and the Markov kernel represents the transition to the new state.

The probability measure $P_{t+1}\in\mathcal{P}(\mathcal{X})$ and the factorial moment measure $M_{t+1}\in\mathcal{M}(\mathcal{X})$ for the Chapman-Kolmogorov equation are found to be, for all $k\in\mathbb{N}$,
\begin{align}
\label{eq:ChapmanKolmogorovProba}
P_{t+1}(dx_{1:k}) & =  \mathbb{E}\left[P_{t+1|t}(dx_{1:k}|\varphi_t)\right],\\
\label{eq:ChapmanKolmogorovMoment}
M_{t+1}(dx_{1:k}) & = \mathbb{E}\left[M_{t+1|t}(dx_{1:k}|\varphi_t)\right].
\end{align}

\subsection{Branching processes}

Branching processes were first proposed by Watson \& Galton \cite{watsongalton} in their theory for natural inheritance to investigate the probability that an aristocratic family name becomes extinct.
The key contribution was to describe the evolution of the distribution of the number of individuals over time as a recursion in terms of its probability generating function (p.g.f.).
Moyal \cite{Moyal62} introduced the probability generating functional (p.g.fl.), which generalises the p.g.f. to account for spatial densities and showed that the
p.g.fl. version of the Galton-Watson functional recursion provides a simple means of describing the Chapman-Kolmogorov equation, where the process is comprised of a number of independent parts, known as multiplicative population processes~\cite{MoyalAP}. The mathematical model is analogous to the model of clustering, where each point in a parent point process generates a new daughter point process.

Let us now consider the transition $\mathbb{E}\left[F(w)(\phi_{t+1})|\phi_t = \{y_1,\ldots,y_m\} \right]$, derived from a general interacting population process (Theorem \ref{thm:hierInteractingPopProc}), but with the following two simplifying assumptions:
\begin{enumerate}
\item \label{it:from1ToN} Transitions $\mathbb{E}\left[ F(w)(\phi_{t+1}) | \omega \right]$ are originated from a single point $\omega = y$ in the prior process $\phi_t$. The p.m.f. $p_t$ is then selected such that $p_t(n) \neq 0$ if and only if $n\leq 1$,
\item \label{it:singleClutterFunc} A single independent process $\mathbb{E}\left[ F(w)(\varphi_{t+1})|\emptyset \right]$ has to be considered so that $P_{\psi}(dz_{1:k})$ is non-zero if and only if $k = m + 1$.
\end{enumerate}

The normalization constant can now be found to be $p_t(0) p_t(1)^m$ so that
\begin{equation*}
\mathbb{E}[F(w)(\phi_{t+1}) | y_{1:m}] = \int\mathbb{E}\left[ F(w)(\varphi_{t+1})|\emptyset,z_{m+1} \right]\prod_{i=1}^m \mathbb{E}\left[ F(w)(\varphi_{t+1}) | y_i,z_i \right] P_{\psi}(dz_{1:m+1}).
\end{equation*}

In its general form, the corresponding transition kernel reads
\begin{equation*}
\bar{P}_{t+1|t}\left( dx_{1:n} | y_{1:m} \right) \propto \sum_{\pi \in \Pi^*(x_{1:n}\cup y_{1:m})} \mathbb{E}_{\psi} \left[Q_P(\pi,\psi)\right],
\end{equation*}
with
\begin{equation*}
Q_P(\pi,z_{1:k}) = \prod_{i=1}^{\dot{\pi}} p_t(\dot{\pi}_{i,x})P_{t+1|t}\left(d\pi_{i,x}|\pi_{i,y},z_i\right)\prod_{i=\dot{\pi}+1}^k p_t(0)P_{t+1|t}\left(\emptyset|\emptyset,z_i\right)
\end{equation*}
and where $\Pi^*(x_{1:n}\cup y_{1:m})$ is the set of all the partition of $\{x_1,\ldots,x_n,y_1,\ldots,y_m\}$ satisfying the following constraints:
\begin{itemize}
\item There is a single set $\{\omega,\emptyset\}$, $\omega \subset \{x_1,\ldots,x_n\}$ (corresponding to assumption \ref{it:singleClutterFunc}) and
\item all the other subsets in the partitions $\pi$ are of one of the two following forms: $\{\omega,y\}$ or $\{\emptyset,y\}$ with $\omega \subset \{x_1,\ldots,x_n\}$ (assumption \ref{it:from1ToN}).
\end{itemize}

Therefore, there are $m+1$ partitions in $\Pi^*(x_{1:n}\cup y_{1:m})$ and since $P_{\psi}(dz_{1:k})$ is non-zero if and only if $k = m + 1$,
\begin{equation*}
\prod_{i=\dot{\pi}+1}^k p_t(0)P_{t+1|t}\left(\emptyset|\emptyset \right) = \prod_{i=m+2}^{m+1} p_t(0)P_{t+1|t}\left(\emptyset|\emptyset \right) = 1.
\end{equation*}

Considering that the normalizing constant is still $p_t(0) p_t(1)^m$, the Markov kernel $\bar{P}_{t+1|t}\left( dx_{1:n} | y_{1:m} \right)$ is
\begin{equation}
\label{eq:MarkovKernelBranching}
\bar{P}_{t+1|t}(dx_{1:n} | y_{1:m}) = \sum_{\pi \in \Pi^*(x_{1:n}\cup y_{1:m})} \int \prod_{i=1}^{m+1} P_{t+1|t}\left(d\pi_{i,x} | \pi_{i,y},z_i\right)P_{\psi}(dz_{1:m+1}).
\end{equation}

Equation (\ref{eq:MarkovKernelBranching}) generalizes branching processes by introducing a weak level of interactions between the different ``branches''. This interaction is weak enough to allow the same kind of simplification as for the usual branching equations, but still allows to, e.g., control the overall number of daughters depending on the number of parents.

\begin{myexample}[Galton-Watson equation]

Branching processes are also associated with the Galton-Watson equation \cite{DaleyVere-Jones,watsongalton} which is a probability generating functional recursion depicting the idea of one parent for any number of daughters and which reads
\begin{equation}
\label{eq:GaltonWatson}
G_{t+1}(h) = G_t(G_{t+1|t}(h|\cdot)).
\end{equation}

The Galton-Watson equation can be recovered from the expression of $\mathbb{E}[F(w)(\phi_{t+1}) |\phi_t = x_{1:k}]$ (Theorem \ref{thm:hierInteractingPopProc}) by considering the parent process $\varphi_t$ and daughter process $\varphi_{t+1}$ to be independent on the process $\psi$ and by considering $p_{\psi}(n) = 1$ if and only if $n = k$ and $p_t(n) = 1$ if and only if $n = 1$ so that
\begin{equation*}
\mathbb{E}[F(w)(\phi_{t+1}) |\phi_t = x_{1:k}] = \prod_{x\in x_{1:k}} \mathbb{E}[F(w)(\varphi_{t+1}) | x].
\end{equation*}

Integrating over the parent process and considering the symmetrised form by taking $F(w) = w^s$, we get the Galton-Watson recursion (\ref{eq:GaltonWatson}).

\end{myexample}

\begin{myexample}[Branching for independent stochastic processes]
\label{ex:BranchingForIndProc}
In this example, a multi-object prior with independent but not identically distributed single-object probability measures is considered, and the conjugacy is demonstrated for this prior following Vo and Vo \cite{VoVoConjugatePrior1,VoVoConjugatePrior2}. The objective here is to derive the same result in a different way, using the generalization of probability generating functionals and the result for branching processes as detailed before.

Let us consider that the probability measure $P_t\in\mathcal{P}(\mathcal{Y})$ of the prior process is defined, for all $m\in\mathbb{N}$, by its projection on $\mathbf{Y}^m$:
\begin{equation}
\label{eq:BayesIndPrior}
P_t(dy_1\times\ldots\times dy_m) = \sum_{\xi \in \Xi(m)} p_{\xi}(m) \prod_{i=1}^m P^{(m)}_{\xi,i}(dy_i),
\end{equation}
where $\Xi(m)$ is a given set of indices and where, for each $\xi \in \Xi(m)$, there exists a parameter set
\begin{equation}
\left\{p_{\xi}(m),\left\{P^{(m)}_{\xi,i}\right\}_{i=1}^m\right\}_{m\geq 0}.
\end{equation}

Since the special case of independent processes is considered here, the Markov kernel (\ref{eq:MarkovKernelBranching}) is simplified to only take into account ``one to one'' or ``one to zero'' transitions which is equivalent to set $P_{t+1|t}(dx_{1:n}|y_{1:m}) \neq 0$, if and only if $n\leq 1$, so that
\begin{equation}
\label{eq:MarkovKernelIndBranching}
\bar{P}_{t+1|t}(dx_{1:n}|y_{1:m}) = \sum_{\pi \in \Pi^*(x\cup y)} P_{t+1|\emptyset}\left(d\nu\right) \prod_{(x_k,y_j)\in\pi} P_{t+1|t}^{(k,j)}\left(dx_k|y_j\right)
\prod_{(\emptyset,y_j)\in\pi} P_{t+1|t}^{(0,j)}\left(\emptyset|y_j\right).
\end{equation}
where $\nu$ is the only subset of $\pi$ such that $\nu\cap\{y_{1:m}\} = \emptyset$ and where $P_{t+1|\emptyset} = P_{t+1|t}(\cdot|\emptyset)$ is the probability measure of the appearing part of the process $\phi_{t+1}$.

Let us assume additionally that the probability measure $P_{t+1|\emptyset}$ is of the form
\begin{equation*}
P_{t+1|\emptyset}(dx_{1:k}) = p_{t+1|\emptyset}(k) \prod_{i=1}^k P^{(k)}_{t+1|\emptyset,i}(dx_i).
\end{equation*}

Writing the Chapman Kolmogorov equation (\ref{eq:ChapmanKolmogorovProba}) with the independent Markov kernel (\ref{eq:MarkovKernelIndBranching}) and the independent prior (\ref{eq:BayesIndPrior}) gives
\begin{equation}
\label{eq:ChapKolIndPosterior}
P_{t+1}(dx_{1:n}) = \sum_{\xi^+ \in \Xi^+(n)} p_{\xi^+}(n) \prod_{i=1}^n P^{(n)}_{\xi^+,i}(dx_i),
\end{equation}
where
\begin{align*}
\Xi^+(n) & = \left\{ \left(m,\xi,\pi\right) \in \left(\mathbb{N}\times\Xi\times\Pi^*(x_{1:n}\cup y_{1:m})\right) \right\},\\
p_{\xi^+}(n) & = p_{\xi}(m)p_{t+1|\emptyset}(\dot{\nu}) \prod_{(x_k,y_j) \in \pi} \int P_{t+1|t}^{(k,j)}\left(d\bar{x} | \bar{y} \right)P^{(m)}_{\xi,j}(d\bar{y})
\prod_{(\emptyset,y_j) \in \pi} \int P_{t+1|t}^{(0,j)}\left(\emptyset | \bar{y} \right)P^{(m)}_{\xi,j}(d\bar{y}),
\end{align*}
and where, denoting by $[\cdot]$ the Iverson bracket, 
\begin{equation*}
P^{(n)}_{\xi^+,i}(dx_i) = \left[\pi_i = (x_k,y_j)\right] \dfrac{\int P_{t+1|t}^{(k,j)}\left(dx_k | \bar{y} \right)P^{(m)}_{\xi,j}(d\bar{y})}{\int P_{t+1|t}^{(k,j)}\left(d\bar{x} | \bar{y} \right)P^{(m)}_{\xi,j}(d\bar{y})}  + \left[\pi_i = (x_k,\emptyset)\right]P^{(\dot{\nu})}_{t+1|\emptyset,k}(dx_k).
\end{equation*}

The posterior (\ref{eq:ChapKolIndPosterior}) is of the same form as the prior (\ref{eq:BayesIndPrior}), proving the conjugacy.
\end{myexample}

\subsection{Renormalization}

Renormalization is used in statistical physics and percolation theory for dealing with correlations of physical systems across different
scales\footnote{
The concept was originally proposed in quantum field theory~\cite{bogolyubov}, where it is still an active area of research.
Whether the approach proposed is applicable in this context will be left for the subject of future work.
}.
In the Ising model, states are assumed to be on a lattice, where vertex $i$ on the lattice has an associated random variable $S_i$ with some ``spin'' value in $\{-1,1\}$.
Kadanoff \cite{kadanoff} suggested that when neighbouring spins are strongly correlated, the mean of a block of spins should not behave very differently than a single spin.
The idea was to calculate physical observables by summing recursively over short-distance degrees of freedom.
Thus, in the Ising model, we can create a  more coarsely grained lattice by only considering nearest-neighbour spins.
This idea was adopted by Wilson \cite{wilson} for considering problems of critical phenomena and phase transitions.

\begin{figure}[htp]
\centering
\includegraphics[width=120px]{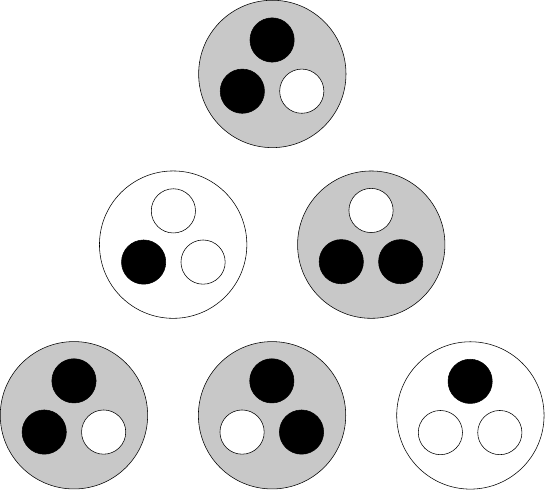}
\vspace{-6pt}
\caption{Renormalization example.}
\label{fig:renorm}
\end{figure}

We describe a simple example of renormalization, adapted from Chowdhury \& Stauffer \cite{chowdhury}.
In Figure \ref{fig:renorm}, the smaller circles represent the original system, and the larger circles represent the renormalised system.
The spin of each smaller circle is denoted by it being either black or white.
The colour of the larger circle is determined by the majority of each of the circles it contains.
Thus, the new system is conditioned on the probabilities of the constituent circles.

\begin{figure}[htp]
\centering
\includegraphics[width=300px]{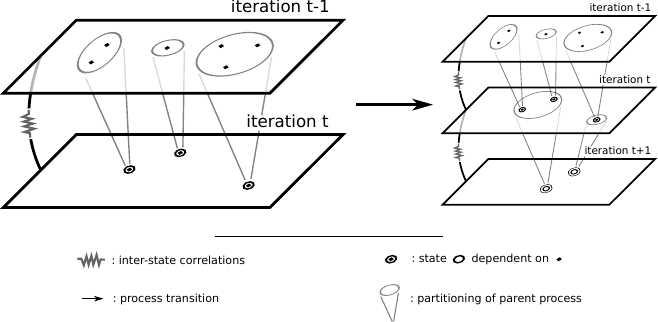}
\caption{Transition for the process corresponding to renormalization.}
\label{fig:renormalization}
\end{figure}

We consider a generalisation of this idea for a system  not restricted to a regular lattice structure.
Let us suppose that we have a system of particles, whose correlations are described with some probability measure $P_x\in\mathcal{P}(\mathcal{X})$.
We then create a new system of  particles, where each new particle is conditioned on groups of parent particles from the original system (see Figure \ref{fig:renormalization}).
This is the converse problem of branching, where groups of particles are conditioned on a single parent process.

\begin{mycorollary}[Renormalization, from Theorem \ref{thm:ProbaAndFactOfHierIntPopProc}]
Let $\phi_{t+1}$ be independent of $\psi$, let the p.m.f. $p_t$ be such that $p_t(0) = 0$ and let the Markov kernel $P_{t+1|t}(dy_{1:m}|x_{1:n})$ be non-zero only if $n=1$, then the conditional expectation of $F(w)(\phi_{t+1})$ given that $\varphi_t = \{x_1,\ldots,x_k\}$ is
\begin{equation}
\label{eq:renormMain}
\mathbb{E}[ F(w)(\phi_{t+1}) | x_{1:k}] \propto \sum_{\pi \in \Pi(x_{1:k})} p_{\psi}(\dot{\pi}) \prod_{\omega\in\pi} p_t(\dot{\omega})\int w(y) P_{t+1|t}(dy|\omega)
\end{equation}
where the normalization constant is found to be
\begin{equation*}
\displaystyle\sum_{\pi \in \Pi(x_{1:k})} p_{\psi}(\dot{\pi}) \prod_{\omega\in\pi} p_t(\dot{\omega}).
\end{equation*}

This transition is a superposition of ``many to one'' transitions. The choice of p.m.f.s $p_{\psi}$ and $p_t$ has still to be done, the former controls the number of independent transitions superposed while the latter controls what are the possible values for the ``many'' in ``many to one'' and their respective probability.
\end{mycorollary}

\begin{myexample}Let us consider the case depicted in Figure \ref{fig:renorm}. Given that there are $k$ particles in the system before renormalization where there exists $r\in\mathbb{N}$ such that $k = 3r$, the different parameters are selected such that the number of transitions is $k/3$ ($p_{\psi}(n) = 1$ only if $n = k/3$) and each transition is composed of $3$ elements of the previous process ($p_t(n) = 1$ only if $n = 3$). It is also assumed that $P_{t+1|t}(dy|\omega)$ is non-zero only if $\omega = \{x_{i_1},x_{i_2},x_{i_3}\}$ where $\{x_{i_1},x_{i_2},x_{i_3}\}$ are neighbours. It results that the single subsisting element in the sum over $\Pi(x_{1:k})$ is the partition of $x_{1:k}$ as in Figure \ref{fig:renorm}. The resulting expectation $\mathbb{E}\left[ F(w)(\phi_{t+1}) | x_{1:k} \right]$ is then
\begin{equation*}
\prod_{\{x_{i_1},x_{i_2},x_{i_3}\} \in x_{1:k}} \int w(y) P_{t+1|t}(dy|x_{i_1},x_{i_2},x_{i_3}).
\end{equation*}

\end{myexample}

\section{Bayesian filtering and smoothing}
\label{sec:BayesRule}

Stochastic filtering is concerned with the on-line estimation of the state of a signal at given intervals, based on a sequence of noisy observations~\cite{Jazwinski}.
Smoothing refines previous filtered estimates based on current measurements.
For discrete-time systems, the stochastic filter, also known as the Bayes' filter, comprises of the Chapman-Kolmogorov equation, for describing the time-evolution of the system, and Bayes' rule, which updates the process based on new measurements.
The corresponding smoother was derived by Kitagawa \cite{kitagawa_original} as a backward recursion.

In this section, we describe Bayesian filtering for the general stochastic population processes introduced in the previous sections.
Finally, we describe the forward-backward smoother for multi-object systems.

\subsection{Bayesian estimation}

Mahler \cite{Mahlermaths} proposed a means of determining the posterior probability generating functional based on a set of measurements.
This idea extends to the general non-symmetric population processes, as described in the following theorem.

\begin{mytheorem}[Bayes' rule]
\label{thm:BayesForGF}
Let $\varphi_x$ be a process with probability measure $P_x\in\mathcal{P}(\mathcal{X})$ and $P_{z|x}$ be a likelihood derived from the Markov kernel $P_{z|x}\in\mathcal{K}(\mathbf{B}_{\mathcal{Z}}\times\mathcal{X})$. The expectation of the updated process is
\begin{equation}
\label{eq:BayesForGF}
\mathbb{E}\left[F(w)(\varphi_x)|z_{1:m}\right] = \dfrac{\mathbb{E}\left[F(w)(\varphi_x) P_{z|x}(z_{1:m}|\varphi_x)\right]}
{\mathbb{E}\left[P_{z|x}(z_{1:m}|\varphi_x)\right]},
\end{equation}
where $F(w)\in\mathcal{B}_b(\mathcal{X})$ is a test function.
\end{mytheorem}

\begin{IEEEproof}
Let $P_{x|z}(\cdot|z_{1:m})\in\mathcal{P}(\mathcal{X})$ be a probability measure, the Bayes' update for $P_{x|z}(\cdot|z_{1:m})$ can be written
\begin{equation*}
P_{x|z}(dx_{1:n}|z_{1:m})\mathbb{E}\left[P_{z|x}(dz_{1:m}|\varphi_x)\right] = P_{z|x}(dz_{1:m}|x_{1:n}) P_x(dx_{1:n}).
\end{equation*}

Assuming the probability measure $P_{z|x}(\cdot|\varphi_x)$ to be absolutely continuous w.r.t. a reference measure $\mu \in\mathcal{M}(\mathcal{Z})$ and using the same notation for the associated density $dP_{z|x}(\cdot|\varphi_x)/d\mu$, the result is obtained when considering the expectation of $F(w)(\varphi_x)$ with probability measure $P_{x|z}(\cdot|z_{1:m})$.
\end{IEEEproof}

\begin{mycorollary}[From Theorem \ref{thm:BayesForGF}]
For all $n\in\mathbb{N}$, the factorial moment measure $M_{x|z}\in\mathcal{M}(\mathcal{X})$ of the updated process is found with
\begin{equation}
\label{eq:cfmdBayes}
M_{x|z}(dx_{1:n}|z_{1:m}) = \dfrac{\mathbb{E}\left[ C_n\left(dx_{1:n}|\varphi_x\right) P_{z|x}(z_{1:m}|\varphi_x) \right]}
{\mathbb{E}\left[P_{z|x}(z_{1:m}|\varphi_x)\right]},
\end{equation}
where $C_n$ is the counting measure from Definition \ref{def:countingMeasure}.
\end{mycorollary}

\begin{myexample}[Generalized PHD filter update] The Probability Hypothesis Density (PHD) filter \cite{Mahlermaths} was developed as a tractable means of propagating a stochastic population process for multiple target tracking applications. A Poisson approximation of the prior population process was made before applying Bayes' rule in order to derive a closed-form recursion in the first-order factorial moment density. However, the Poisson approximation ensures that the correlations between the objects are not maintained. We can extend this idea here using a generalized Poisson prior, or Khinchin prior, to obtain a generalized PHD filter update that maintains target correlations.

Let us consider again a Khinchin process as in Section \ref{ssec:SuperProcAndKhinchinProc} and more specifically a joint Khinchin process as in Example \ref{ex:JointKhinchin}. The $k^{th}$-order factorial moment measure of an updated Khinchin process via Bayes' rule can be written
\begin{equation}
\label{KhinchinProcBayesUpdate}
M(dx_{1:k}|z_{1:m}) = \dfrac
{\displaystyle\sum_{\pi\in\Pi(x\cup z)} \prod_{\omega\in\pi} \sum_{n\geq 0}
\int \sideset{}{^{\neq}}\sum_{1\leq j_1,\ldots,j_{\dot{\omega}_x} \leq n} J(\omega_z|\hat{x}^{j_{1:\dot{\omega}_x}}_{1:n})\lambda_n(d\hat{x}^{j_{1:\dot{\omega}_x}}_{1:n})}
{\displaystyle\sum_{\pi\in \Pi(z)} \prod_{\omega \in \pi} \left( \nu(\omega) + \sum_{n\geq1} \int J(\omega|\bar{x}_{1:n})\lambda_n(d\bar{x}_{1:n}) \right)},
\end{equation}
where $J(\omega_z|\emptyset)\lambda_0 = \nu(\omega_z)$ represents spurious measurements, and for all $r\in\mathbb{N}$, $1\leq r \leq n$, $\displaystyle\hat{x}^{j_{1:\dot{\omega}_x}}_r$ is defined to be
\begin{align}
\left\{
\begin{array}{ll}
\omega_{x,i} & \text{ if } \exists i\in\mathbb{N}, 1\leq i\leq\dot{\omega}_x, r = j_i \\
\bar{x}_r & \text{ otherwise},
\end{array}
\right.
\end{align}
and where the denominator is obtained by differentiation of (\ref{eq:JointKhinchin}) w.r.t. to the second argument at point $0$ and considering the first argument at point $1$.

This expression is very general and is difficult to use as it is, but simpler cases can be deduced straightforwardly from it as in the following examples concerning the PHD filter and a filter handling occlusions and unresolved objects.
\end{myexample}

\begin{myexample}[The PHD filter update]
Starting from a the Bayes' update of a Khinchin process (\ref{KhinchinProcBayesUpdate}) written at the first order ($k=1$), the assumptions needed to arrive to the usual PHD filter are that
\begin{enumerate}
\item Poisson measures $\lambda_n$ are zero for $n>1$,
\item no more than one measurement is originated from each object and
\item measurements are independent of each other.
\end{enumerate}

The only partitions left in the set $\Pi(z_{1:m}\cup x_{1:k})$ are of the following forms:
\begin{itemize}
\item $\{x,z_1,\ldots,z_m\}$,
\item $\{z_1,\ldots,\{x,z_i\},\ldots,z_m\}$, $1\leq i \leq m$.
\end{itemize}

Therefore, noting that $J(\omega_z|x) = P(\omega_z|x)$ when $\dot{\omega}_z \leq 1$, (\ref{KhinchinProcBayesUpdate}) becomes the PHD filter update:
\begin{equation*}
M(dx|z_{1:m}) = P(\emptyset|x)\lambda_1(dx) + \sum_{z\in z_{1:m}} \dfrac{P(z|x)\lambda_1(dx)}{\nu(z) + \int P(z|\bar{x})\lambda_1(d\bar{x})}.
\end{equation*}
\end{myexample}

To find other kind of filters that are likely to be tractable in practice, one can look for assumptions that are going to simplify (\ref{KhinchinProcBayesUpdate}) without being as strong as for the PHD filter case. It is for instance possible to only assume that measurements are independent of each other so that the sum over partitions in the denominator of (\ref{KhinchinProcBayesUpdate}) reduces to one term. This assumption is studied in the next example.

\begin{myexample}[Filter for occlusions and unresolved objects]
It is assumed here that each joint states in the Khinchin process can generate no more than one measurement. From a physical point of view, this assumption means that some group of object might only generate one measurement because of occlusions or because these objects are close of each other when compared to the sensor resolution (unresolved objects).

In this case, (\ref{KhinchinProcBayesUpdate}) simplifies to
\begin{equation*}
M(dx_{1:k}|z_{1:m}) = \sum_{\pi\in\Pi^*(x\cup z)}\dfrac
{\displaystyle\prod_{\omega\in\pi} \sum_{n\geq 0} \int \sideset{}{^{\neq}}\sum_{1\leq j_1,\ldots,j_{\dot{\omega}_x} \leq n} P(\omega_z|\hat{x}^{j_{1:\dot{\omega}_x}}_{1:n})\lambda_n(d\hat{x}^{j_{1:\dot{\omega}_x}}_{1:n})}
{\displaystyle\prod_{z \in z_{1:m}} \left( \nu(z) + \sum_{n\geq1} \int P(z|\bar{x}_{1:n})\lambda_n(d\bar{x}_{1:n}) \right)},
\end{equation*}
where subsets of the partition $\pi\in\Pi^*(x\cup z)$ have no more than one element in the individual measurement space $\mathbf{Z}$ and where $\hat{x}^{j_{1:\dot{\omega}_x}}_{1:n}$ is defined as before.

Let us now reduce the Khinchin prior to a Gauss-Poisson prior, for which $\lambda_n$ is non-zero only if $n\leq 2$. For the first moment ($k=1$), the only partitions left are
\begin{itemize}
\item $\{x,z_1,\ldots,z_m\}$,
\item $\{z_1,\ldots,\{x,z_i\},\ldots,z_m\}$, $1\leq i \leq m$.
\end{itemize}
as for the PHD filter.

The first-order moment measure $M(dx|z_{1:m})$ can be written
\begin{equation*}
M(dx|z_{1:m}) = F(\emptyset,x) + \sum_{z\in z_{1:m}}\dfrac{F(z,x)}{C(z)}
\end{equation*}
where, for $\omega \in \{\emptyset,z_1,\ldots,z_m\}$,
\begin{equation*}
F(\omega,x) = P(\omega|x)\lambda_1(dx) + 2\int P(\omega|x,\bar{x})\lambda^s_2(dx\times d\bar{x}),
\end{equation*}
with $2\lambda^s_2(dx_1\times dx_2) = \sum_{\sigma}\lambda_2(dx_{\sigma_1}\times dx_{\sigma_2})$ and
\begin{equation*}
C(z) = \nu(z) + \int P(z|\bar{x})\lambda_1(d\bar{x}) + \int P(z|\bar{x}_1,\bar{x}_2)\lambda^s_2(d\bar{x}_1\times d\bar{x}_2).
\end{equation*}

We can also determine the second-order factorial moment measure which depicts correlations between objects:
\begin{multline*}
M(dx_1\times dx_2|z_{1:m}) = F_2(\emptyset,x_1,x_2) + F(\emptyset,x_1)F(\emptyset,x_2) \\
+ \sum_{z\in z_{1:m}} \left(F(\emptyset,x_1)\dfrac{F(z,x_2)}{C(z)} + F(\emptyset,x_2)\dfrac{F(z,x_1)}{C(z)}\right) + \sum_{z\in z_{1:m}} \dfrac{F_2(z,x_1,x_2)}{C(z)} + \sideset{}{^{\neq}}\sum_{z,z' \in z_{1:m}} \dfrac{F(z,x_1)F(z',x_2)}{C(z)C(z')},
\end{multline*}
where $F_2(\omega,x_1,x_2) = 2P(\omega|x_1,x_2)\lambda^s_2(dx_1\times dx_2)$.

An alternative presentation of filters with a Gauss-Poisson prior can be found in \cite{SumeetFiltersForSpatialPP}.
\end{myexample}

\begin{myexample}[Independent stochastic processes]
In this example, the conjugacy of independent priors through Bayes' rule is demonstrated in a way similar to  Example \ref{ex:BranchingForIndProc}.

Rewriting the Markov kernel (\ref{eq:MarkovKernelBranching}) with adapted notations for Bayes' rule and assuming the independence w.r.t. the process $\psi$, this kernel becomes
\begin{equation}
\label{eq:MarkovKernelBayesRule}
\bar{P}_{z|x}\left( dz_{1:m} | x_{1:n} \right) = \sum_{\pi \in \Pi^*(z\cup x)} P_z(d\nu)\prod_{(\omega_z,x_k)\in\pi} P^{(k)}_{z|x}\left(d\omega_z | x_k \right),
\end{equation}
where $\Pi^*(z\cup x)$ is the set of all partitions of $\{z_{1:m},x_{1:n}\}$ satisfying the following constraints: there is a single set $\{\nu,\emptyset\}$, $\nu \subset x_{1:n}$ and all the subsets in the partitions $\pi$ take the forms $\{\omega,x\}$ or $\{\emptyset,x\}$ with $\omega \subset z_{1:m}$. $P_z = P_{z|\emptyset}$ is the probability measure of spurious measurements.

It is assumed that $\bar{P}_{z|x}$ and $P^{(k)}_{z|x}$ are absolutely continuous w.r.t. a reference measure $\mu\in\mathcal{M}(\mathcal{Z})$ and the same notation are used for the associated densities. Applying a $n^{th}$-order partial differentiation on (\ref{eq:BayesForGF}) w.r.t. $\{\mathbf{1}_{dx_1},\ldots,\mathbf{1}_{dx_n}\}$ at point $w=0$ and using the Markov kernel (\ref{eq:MarkovKernelBayesRule}), Bayes' rule reads
\begin{equation*}
P^+(dx_{1:n}|z_{1:m}) = C^{-1}(z_{1:m})\sum_{\pi \in \Pi^*(z\cup x)} P_z(\nu)\prod_{(\omega_z,x_k)\in\pi} P^{(k)}_{z|x}\left(\omega_z | x_k \right)P(dx_{1:n}),
\end{equation*}
where $C(z_{1:m})$ is defined to be
\begin{equation*}
\sum_{k\geq 0}\int \sum_{\pi \in \Pi^*(z\cup \bar{x})} P_z(\nu)\prod_{(\omega_z,x_i)\in\pi} P^{(i)}_{z|x}\left(\omega_z | \bar{x}_i \right)P(d\bar{x}_{1:k}).
\end{equation*}

Considering that the prior $P$ is of the same form as (\ref{eq:BayesIndPrior}), the posterior measure $P^+$ can be written
\begin{equation}
\label{eq:BayesIndPosterior}
P^+(dx_{1:n}|z_{1:m}) = \sum_{\xi^+ \in \Xi^+(n)} p_{\xi^+}(n) \prod_{i=1}^n P^{(n)}_{\xi^+,i}(dx_i)
\end{equation}
where
\begin{align*}
\Xi^+(n) & = \left\{ \left(\xi,\pi\right) \in \left(\Xi\times\Pi^*(z_{1:m}\cup x_{1:n})\right) \right\},\\
p_{\xi^+}(n) & = \dfrac{p_{\xi}(n) P_z(\nu)}{C(z_{1:m})} \prod_{\substack{(\omega_z,x_k)\in\pi\\\omega_z\neq\emptyset}} \int P^{(k)}_{z|x}\left(\omega_z | \bar{x} \right)P^{(n)}_{\xi,k}(d\bar{x}),\\
P^{(n)}_{\xi^+,i}(dx_i) & = \left[\pi_i=(\emptyset,x_k)\right] P^{(k)}_{z|x}\left(\emptyset | x_k \right)P_{\xi,k}(dx_k) \\
&\qquad\qquad\qquad\qquad+ \left[\pi_i = (\omega_z,x_k),\omega_z\neq\emptyset\right] \dfrac{P_{z|x}^{(k)} \left(\omega_z | x_k \right)P^{(n)}_{\xi,k}(dx_k)}{\int P_{z|x}^{(k)}\left(\omega_z | \bar{x} \right)P^{(n)}_{\xi,k}(d\bar{x})}.
\end{align*}

The posterior (\ref{eq:BayesIndPosterior}) is the probability measure of an independent process as the prior is, proving the conjugacy.
\end{myexample}

\vskip10pt

Other examples of application can be easily derived. For instance,
higher order interacting stochastic population processes can be defined in a very similar way to handle simultaneous multi-object estimation and multi-sensor registration \cite{RisticCalibration,MahlerRegistration}.

\subsection{Forward-backward smoothing}

The forward-backward smoother, proposed by Kitagawa \cite{kitagawa_original}, refines previous filtered estimates based on measurements up to the current time-step.
This has been studied for multi-object systems for multi-target tracking \cite{MahlerVoVo,fusionsmooth,Hernandez,VoVoMahler}.

\begin{mytheorem}
\label{thm:BayesForGFfb}
Let $\varphi_{t\,\uparrow\, t}$, $\varphi_{t'\,\uparrow\, t'}$ and $\varphi_{t\,\uparrow\, t'}$ be processes with respective measures $P_{t\,\uparrow\, t}\in\mathcal{P}\left(\mathcal{X}_t\right)$, $P_{t'\,\uparrow\, t'} \in \mathcal{P}\left(\mathcal{X}_{t'}\right)$ and $P_{t\,\uparrow\, t'}\in \mathcal{P}\left(\mathcal{X}_{t}\right)$, where ``$t_1\,\uparrow\, t_2$'' means ``at time $t_1$ given measurements up to time $t_2$''. Assuming $t<t'$ and $F(w)\in\mathcal{B}_b(\mathcal{X}_t)$, the expectation of the smoothed process is
\begin{equation*}
\mathbb{E}\left[F(w)(\varphi_{t\,\uparrow\, t'})\right] = \mathbb{E}\left[ \mathbb{E}\left[ F(w)(\varphi_{t\,\uparrow\, t'}) | \varphi_{t'\,\uparrow\, t'} \right] \right],
\end{equation*}
where, for $x_i\in\mathbf{X}_{t'}, 1\leq i \leq k$,
\begin{equation*}
\mathbb{E}\left[ F(w)(\varphi_{t\,\uparrow\, t'}) | x_{1:k} \right] = \dfrac{\mathbb{E}\left[F(w)(\varphi_{t\,\uparrow\, t})\, P_{t'|t}\left(x_{1:k}|\varphi_{t\,\uparrow\, t}\right)\right]}
{\mathbb{E}\left[P_{t'|t}(x_{1:k}|\varphi_{t\,\uparrow\, t})\right]},
\end{equation*}
where $P_{t'|t}$ is derived from the kernel $P_{t'|t}\in\mathcal{K}(\mathbf{B}_{\mathcal{X}_{t'}}\times\mathcal{X}_t)$.
\end{mytheorem}

\begin{IEEEproof}
The forward-backward smoother can be viewed as an application of both the Chapman-Kolmogorov equation and Bayes' rule~\cite{fusionsmooth,slamsmoother}.
In particular, the posterior probability measure $P_{t\,\uparrow\, t'}\in\mathcal{P}(\mathcal{X}_t)$, $t<t'$, is found with the Chapman-Kolmogorov equation, where the backward transition is determined by using Bayes' rule with the forward transition as a likelihood:
\begin{align*}
P_{t\,\uparrow\, t'}\left(dx_{1:m}\right) = \mathbb{E}\left[
 \dfrac{P_{t'|t}(\varphi_{t'\,\uparrow\, t'}|x_{1:m} )P_{t\,\uparrow\, t}\left(dx_{1:m}\right)}
 {\mathbb{E}\left[P_{t'|t}(\varphi_{t'\,\uparrow\, t'}|\varphi_{t\,\uparrow\, t})\right]}
\right].
\end{align*}

The result is the expectation of $F(w)(\varphi_{t\,\uparrow\, t'})$.
\end{IEEEproof}

\section{Conclusion}
\label{sec:conclusion}

This paper develops a unified general framework for the estimation of multi-object dynamical systems with interactions and hierarchies.
The generality of the approach is highlighted through the application to a number of case studies.
Specific contributions include
\begin{enumerate*}[label=(\roman{*})]
\item Fa\`a di Bruno's formula for G\^ateaux differentials,
\item the generalisation of probability generating functional from point process theory to account for non-symmetric measures,
\item the unification of the theory branching processes and renormalization in statistical physics and percolation theory through the introduction of interacting population processes, and
\item a general solution to the problem of multi-object estimation with applications to multi-target tracking.
\end{enumerate*}

\section*{Acknowledgements}

Daniel Clark is a Royal Academy of Engineering/EPSRC Research Fellow. Jeremie Houssineau has a PhD scholarship sponsored by DCNS and a tuition fees scholarship by Heriot-Watt University. This work was supported by the Engineering and Physical Sciences Research Council grant EP/H010866/1.

\bibliography{bibHierSys}

\end{document}